\documentclass[11pt]{article}

\makeatletter
\@addtoreset{equation}{section}
\makeatother

\input epsf

\usepackage{amssymb}
\usepackage{amsmath}
\usepackage{epsfig}

\def\tfrac#1#2{{{\lower.6ex
\hbox{$\scriptstyle#1$}}\over
{\raise.7ex
\hbox{$\scriptstyle#2$}}}}

\def\bigO{{\cal O}}
\def\NN{{\mathbb N}}

\def\ZZ{{\mathbb Z}}
\def\calC{{\cal C}}
\def\calL{{\cal L}}
\def\wt{\widetilde}
\def\phase{{\rm ph}}
\def\dsp#1{\displaystyle#1}

\def\protectbold#1{\protect{\boldmath{$#1$}}}

\begin{document}

\title{Large Degree Asymptotics of Generalized Bernoulli and Euler Polynomials}

\author{Jos\'e Luis L\'opez \\
    Departamento de Mat\'ematica e Inform\'atica,\\
    Universidad P\'ublica de Navarra, 31006 Pamplona, Spain\\
       \and
    Nico M. Temme\\
    CWI, Science Park 123, 1098 GX Amsterdam, The Netherlands. \\
     { \small e-mail: {\tt
    jl.lopez@unavarra.es, 
    Nico.Temme@cwi.nl}}
    }

\date{\today}
\maketitle

\begin{abstract}
\noindent
Asymptotic expansions are given for large values of $n$ of the generalized Bernoulli polynomials $B_n^\mu(z)$ and Euler polynomials $E_n^\mu(z)$. In a previous paper L\'opez and Temme (1999) these polynomials have been considered for large values of $\mu$, with $n$ fixed. In the literature no complete description of the large $n$ asymptotics of the considered polynomials is available. We give the general expansions, summarize known results of special cases and give more details about  these results. We use two-point Taylor expansions for obtaining new type of expansions. The analysis is based on contour integrals that follow from the generating functions of the polynomials.
\end{abstract}

\vskip 0.8cm \noindent
{\small
2000 Mathematics Subject Classification:
11B68, 30E10, 33E20, 41A60.
\par\noindent
Keywords \& Phrases:
asymptotic expansions, 
generalized Bernoulli polynomials,
generalized Euler polynomials.

}
\section{Introduction}\label{sec:intro}

Generalized Bernoulli and Euler polynomials of degree
$n$, complex order $\mu$ and complex argument $z$, denoted
respectively by $B_n^{\mu}(z)$ and $E_n^{\mu}(z)$, can be defined by their generating functions. We have \cite{Milne:1951:CFD,Temme:1996:SFI}
\begin{equation}\label{Berdef}
\frac{w^{\mu}e^{wz}}{ (e^w-1)^{\mu}}=
\sum_{n=0}^\infty \frac {B_n^{\mu}(z)}{ n!}w^n,
\quad \vert w\vert<2\pi,
\end{equation}
and
\begin{equation}\label{Euldef}
\frac{2^{\mu}e^{wz}}{ (e^w+1)^{\mu}}=
\sum_{n=0}^\infty \frac{E_n^{\mu}(z)}{ n!}w^n,
\quad \vert w\vert<\pi.
\end{equation}

These polynomials
play an important role in the calculus of finite
differences. In fact,
the coefficients  in all the usual central-difference formulae 
for interpolation,
numerical differentiation and integration, and differences in 
terms of derivatives
can be expressed in terms of these polynomials (see \cite{Milne:1951:CFD,Norlund:1924:VUD}).

An explicit formula for the generalized Bernoulli polynomials can be
found in \cite{Srivastava:1988:EFB}. Properties and explicit formulas for the generalized
Bernoulli and Euler numbers can be found in
\cite{Lundell:1987:OTD,Todorov:1985:UFS,Todorov:1993:ERF} and in cited references.

In a previous paper \cite{Lopez:1999:HPA} we have considered these polynomials for large values of $\mu$, with $n$ fixed, and in the present paper we consider $n$ as the large parameter, with the other parameters fixed. We summarize known results from the literature for integer values of $\mu$, and give more details about these 
results. We describe the method for obtaining the coefficients in the expansion for general $\mu$. Finally, we use two-point Taylor expansions for obtaining new type 
of expansions  for general $\mu$. 
The analysis is based on contour integrals that follow from the generating functions of the polynomials.

\section{The generalized Bernoulli polynomials}\label{sec:Ber}
Three different cases arise, depending on $\mu=0,-1,-2,\ldots$, $\mu=1,2,3,\dots$, and $\mu$ otherwise, real or complex. Our approach is based on the Cauchy integral
\begin{equation}\label{Berint}
B_n^{\mu}(z)=\frac{n!}{ 2\pi i} \int_{\calC}{w^{\mu}e^{wz}\over(e^w-1)^{\mu}}
\frac{dw}{ w^{n+1}},
\end{equation}
where $\calC$ is a circle around the origin, with radius less than $2\pi$. This follows from \eqref{Berdef}.

\subsection{Asymptotic form when  \protectbold{\mu=0,-1,-2,\ldots}}\label{subsec:Ber1}
In this  case the generating series in \eqref{Berdef} converges for all finite values of $z$, and, hence, the polynomials  $B_n^{\mu}(z)$ have a completely different behavior compared with the general case.

We first follow the approach given in  \cite{Weinmann:1963:AEB}, and observe that when $\mu$ is a negative integer or zero, say $\mu=-m$ $(m=0,1,2,\ldots)$, we can express $B_n^{\mu}(z)$ in terms of a finite sum. We expand by the binomial theorem
\begin{equation}\label{binom}
\left(e^w-1\right)^m=\sum_{r=0}^m(-1)^{m-r}{m\choose r}e^{rw}.
\end{equation}
This gives
\begin{equation}\label{Bersum}
B_n^{-m}(z)=\frac{n!}{(n+m)!} \sum_{r=0}^m (-1)^{m-r}{m\choose r} (z+r)^{n+m}.
\end{equation}
For any given $m\in \NN$ and complex $z$ only the term or terms with the largest values of $|z+r|$ will give a large contribution to the sum in \eqref{Bersum}, the other terms being exponentially small in comparison. We conclude, that \eqref{Bersum} gives the asymptotic form when $n\to\infty$, when $\mu=-m$ and $z$ are fixed. 

In particular, when $z>0$, the term with index $r=m$ is maximal, and we have
\begin{equation}\label{Bersum0}
B_n^{-m}(z)=\frac{n!}{(n+m)!}   (z+m)^{n+m}\left[1+\bigO\left(\frac{z+m-1}{z+m}\right)^{n+m}\right].
\end{equation}
The error term can also be estimated by $\bigO(\exp(-(n+m)/(z+m)))$, which is indeed exponentially small compared with unity.

For general complex $z=x+iy$ and $x>-m/2$ the term with index $r=m$ again is maximal and the same estimate as in \eqref{Bersum0} is valid. When $x=-m/2$ the terms with $r=0$ and $r=m$ give the maximal contributions, and we have
 \begin{equation}\label{Bersum1}
B_n^{-m}(z)\sim\frac{n!}{(n+m)!}   \left[(-1)^m(-\tfrac12m+iy)^{n+m}+(\tfrac12m+iy)^{n+m} \right].
\end{equation}
When $x<-m/2$ the term with index $r=0$ is maximal, and we have
\begin{equation}\label{Bersum2}
B_n^{-m}(z)=(-1)^m\frac{n!}{(n+m)!}   z^{n+m}\left[1+\bigO\left(\frac{z+1}{z}\right)^{n+m}\right].
\end{equation}

By using the saddle point method we can obtain an estimate similar as the one in \eqref{Bersum0}. We write
\begin{equation}\label{Bersum3}
B_n^{-m}(z)=\frac{n!}{2\pi i}\int_{{\cal C}}\left(1-e^{-w}\right)^m e^{\phi(w)}\frac{dw}{w}, 
\end{equation}
where $\phi(w)=(z+m)w-(n+m)\ln w$. This function has a saddle point at $w_0=(n+m)/(z+m)$, and when $\Re w_0$ is large and positive we replace $(1-e^{-w})$ by its value at this point. Then we have
\renewcommand{\arraystretch}{1.75}
\begin{equation}\label{Bersum4}
\begin{array}{lll}
B_n^{-m}(z)&\sim&\dsp{\left(1-e^{-w_0}\right)^m\frac{n!}{2\pi i}\int_{{\cal C}} e^{(z+m)w}\frac{dw}{w^{n+m+1}}}\\
&=&
\dsp{\left(1-e^{-w_0}\right)^m \frac{n!}{(n+m)!}(z+m)^{n+m},}
\end{array}
\renewcommand{\arraystretch}{1.0}
\end{equation}
in which $e^{-w_0}$ is exponentially small. When  $\Re w_0$ is not positive we can modify this method to obtain the estimates as given in  \eqref{Bersum1} and \eqref{Bersum2}. Also, we can make further steps in the saddle point analysis, and show that the next terms in the expansion are exponentially small compared with unity, similar as  shown in \eqref{Bersum0}. However, the representation in \eqref{Bersum} describes very elegantly the asymptotic behavior.

\subsection{Asymptotic form when  \protectbold{\mu= 1, 2,3,\ldots}}\label{subsec:Ber2}
We write $\mu=m$. The starting point is the expansion for $m=1$:
\begin{equation}\label{Ber1ser}
B_n^{1}(z)=B_n(z)=-n!\sum_{{k=-\infty \atop k\ne0}}^\infty\frac{e^{2\pi ikz}}{(2\pi ik)^n},
\end{equation}
which for $z\in(0,1)$ can be viewed as a Fourier expansion of $B_n^{1}(z)$. This expansion follows from  taking the radius of the circle $\calC$ in \eqref{Berint} equal to $(2K+1)\pi$, ($K$ an integer). Taking $K$ large, we take into account the poles of the integrand at $w=2\pi ik$ ($k=\pm1,\pm2,\ldots$), and calculate the residues of these poles. The integral around the circle $\calC$ tends to zero as $K\to\infty$, provided $n>1$ and $0\le z\le 1$.

This gives the expansions ($n=1,2,3,\ldots; 0\le z\le1$)
\begin{equation}\label{Bersereven}
B_{2n}(z)=2(-1)^{n+1}(2n)!\,\sum_{k=1}^\infty\frac{\cos(2\pi kz)}{(2\pi k)^{2n}},
\end{equation}
and
\begin{equation}\label{Berserodd}
B_{2n+1}(z)=2(-1)^{n+1}(2n+1)!\,\sum_{k=1}^\infty\frac{\sin(2\pi kz)}{(2\pi k)^{2n+1}}.
\end{equation}
In \eqref{Berserodd} we can take $n=0$, provided $0<z<1$. This gives the well-known Fourier expansion of  $B_{1}(z)=z-1/2$, $0<z<1$.

In \eqref{Ber1ser} only the terms  with $k=\pm1$ are relevant for the asymptotic behavior, and we obtain for fixed complex $z$
\begin{equation}\label{Ber1as}
B_n^{1}(z)=\frac{2(-1)^{n+1}n!}{(2\pi)^n}\left[\cos\left(2\pi z+\tfrac12\pi n\right) +\bigO\left(2^{-n}\right)\right],\quad n\to\infty.
\end{equation}

For general fixed real or complex $z$ the series in  \eqref{Bersereven} and  \eqref{Berserodd} can be viewed as asymptotic expansion for large $n$, as easily follows from the ratio test.

For general $\mu=m=1,2,3,\ldots$ a similar expansion as in \eqref{Ber1ser} can be given. In that case the poles of the integrand in  \eqref{Berint} are of higher order. We can write
\begin{equation}\label{Bermser}
B_n^{m}(z)=-n!\sum_{{k=-\infty \atop k\ne0}}^\infty \beta_k^m(n,z)\frac{e^{2\pi ikz}}{(2\pi ik)^n},
\end{equation}
where $\beta_k^1(n,z)=1, \forall k$.  This is a Fourier expansion for $z\in(0,1)$ when $m<n$. For other values of $z$ it can be used as an asymptotic expansion for large $n$.

An explicit form of   $\beta_k^{m}(n,z)$ follows from calculating  the residues of the poles at $2\pi ik$ of order $m$ of the integrand in \eqref{Berint}. For this we compute the coefficient $c_{m-1}$ in the expansion
\begin{equation}\label{Berres1}
\frac{(w-2\pi ik)^m w^m e^{zw}}{(e^w-1)^mw^{n+1}}=\sum_{r=0}^\infty c_r (w-2\pi ik)^r,
\end{equation}
from which we obtain
\begin{equation}\label{Bermser1a}
\beta_k^m(n,z)\frac{e^{2\pi ikz}}{(2\pi ik)^n}=c_{m-1}.
\end{equation}
We substitute  $w=s+2\pi ik$ and write the expansion as
\begin{equation}\label{Berres2}
e^{2\pi ikz}\frac{s^m e^{zs}}{(e^s-1)^m(s+2\pi ik)^{n+1-m}}=\sum_{r=0}^\infty c_r s^r.
\end{equation}
We use \eqref{Berdef} and write the left-hand side in the form
\begin{equation}\label{Berres2a}
\frac{e^{2\pi ikz}}{(2\pi ik)^{n+1-m}}\sum_{\nu=0}^\infty \frac {B_\nu^{m}(z)}{ n!}s^\nu 
\sum_{\nu=0}^\infty{m-n-1\choose \nu}\frac{s^\nu}{(2\pi ik)^\nu}.
\end{equation}
Hence, $c_r$ of \eqref{Berres2} can be written as
\begin{equation}\label{Berres2b}
c_{r}=\frac{e^{2\pi ikz}}{(2\pi ik)^{n+1-m}}\sum_{\nu=0}^{r} \frac{B_\nu^m(z)}{\nu!}{m-n-1\choose r-\nu}(2\pi ik)^{\nu-r},
\end{equation}
and we conclude that 
\begin{equation}\label{Berres3}
c_{m-1}=\frac{e^{2\pi ikz}}{(2\pi ik)^{n}}\sum_{\nu=0}^{m-1} \frac{B_\nu^m(z)}{\nu!}{m-n-1\choose m-1-\nu}(2\pi ik)^{\nu}.
\end{equation}
It follows from \eqref{Bermser1a}  that 
\begin{equation}\label{Berres4}
\beta_k^m(n,z)=\sum_{\nu=0}^{m-1} \frac{B_\nu^m(z)}{\nu!}
{m-n-1\choose m-1-\nu}(2\pi ik)^{\nu}.
\end{equation}
To avoid binomials with negative integers, and to extract the main asymptotic factor, we write
\begin{equation}\label{Berres5}
\beta_k^m(n,z)=(-1)^{m-1}{n-1\choose m-1}\sum_{\nu=0}^{m-1} B_\nu^m(z)
{m-1\choose \nu}\frac{(n-\nu-1)!}{(n-1)!}(-2\pi ik)^{\nu}.
\end{equation}

For large $n$ the main term occurs for $\nu=0$. We have
\begin{equation}\label{Berres6}
\beta_k^m(n,z)=\frac{(-1)^{m-1}n^{m-1}}{(m-1)!}\left[1+\bigO(n^{-1})\right].
\end{equation}

Observing that, as in \eqref{Ber1ser}, only the terms with $k=\pm1$ are relevant for the asymptotic behavior, we obtain 
\begin{equation}\label{Bermser2}
B_n^{m}(z)=\frac{(-1)^{n+1}n!}{(2\pi)^n}\left[\beta_1^{m}(n,z)e^{2\pi iz+\frac12\pi in}+
\beta_{-1}^{m}(n,z)e^{-2\pi iz-\frac12\pi in}+\ldots\right],
\end{equation}
and by using \eqref{Berres5} we obtain for fixed $m$ and complex $z$ (cf. \eqref{Ber1as})
\renewcommand{\arraystretch}{1.5}
\begin{equation}\label{Bermas1}
\begin{array}{ll}
\dsp{B_n^{m}(z)=\frac{2(-1)^{m+n}}{(2\pi)^n}{n-1\choose m-1} \,\times} \\
\quad\quad\quad
\dsp{\left[
\sum_{\nu=0}^{m-1} B_\nu^m(z)
{m-1\choose \nu}\frac{(n-\nu-1)!}{(n-1)!}(2\pi )^{\nu}\cos\sigma+\bigO\left(2^{-n}\right)\right]},
\end{array}
\renewcommand{\arraystretch}{1.0}
\end{equation}
as $n\to\infty$, where $\sigma=(2z+\tfrac12n-\tfrac12\nu)\pi$.

To obtain $\beta_k^m(n,z)$ for  $m> 1$ we can also use a recurrence relation. We have the relation
\begin{equation}\label{Berrec}
\mu B_n^{\mu+1}(z)=(\mu-n)B_n^{\mu}(z)+n (z-\mu)B_{n-1}^{\mu}(z),\quad n\ge 1,
\end{equation}
which follows from \eqref{Berdef} by differentiating both members with respect to $w$. By differentiation with respect to $z$ we find
\begin{equation}\label{Berdif}
 nB_{n-1}^{\mu}(z)=\frac{d}{dz}B_n^{\mu}(z),
\end{equation}
giving
\begin{equation}\label{Berrecd}
\mu B_n^{\mu+1}(z)=(\mu-n)B_n^{\mu}(z)+ (z-\mu)\frac{d}{dz}B_n^{\mu}(z),\quad n\ge 0.
\end{equation}
This gives the recurrence relation for $m=1,2,3,\ldots$
\begin{equation}\label{betarec}
m \beta_k^{m+1}(n,z)=[m-n+2\pi ik(z-m)]\beta_k^{m}(n,z)+ (z-m)\frac{d}{dz}\beta_k^{m}(n,z).
\end{equation}

\subsection{Asymptotic form for general complex  \protectbold{\mu}}\label{subsec:Ber3}
\begin{figure}
\caption{\small Contour for \eqref{Berint} for general $\mu$
\label{Ber.fig1}}
\begin{center}
\epsfxsize=8cm \epsfbox{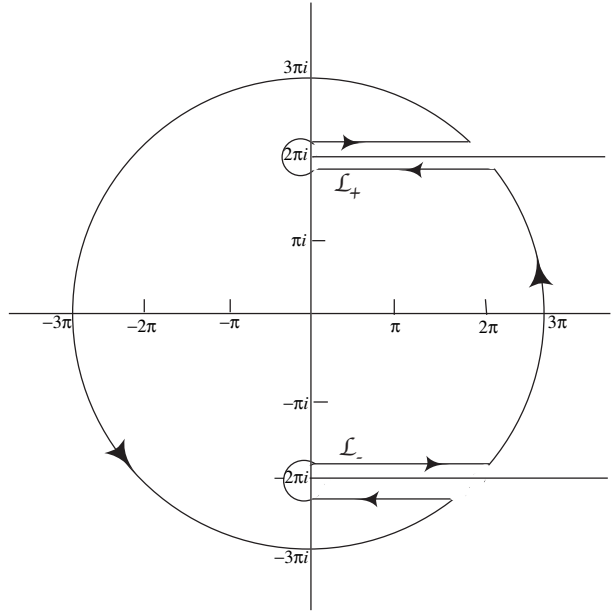}
\end{center}
\end{figure}

We consider \eqref{Berint} and observe that the singularities at $\pm2\pi i$ are the sources for  the main asymptotic contributions. We integrate around a circle with radius $3\pi$, avoiding branch cuts running from $\pm2\pi i$ to $+\infty$. See Figure~\ref{Ber.fig1}. The contribution from the circular arc is $\bigO((3\pi)^{-n})$, which is exponentially small with respect to the main contributions.

We denote the loops by $\calL_{\pm}$ and the contributions from the loops by $I_{\pm}$. For the upper loop  we substitute $w=2\pi i e^s$. This gives
\begin{equation}\label{Ip}
I_+=\frac{n!}{2\pi i}\frac{e^{2\pi iz}}{(2\pi i)^n}\int_{\calC_+} g(s) s^{-\mu}e^{-ns}\,ds,
\end{equation}
where
\begin{equation}\label{gs}
g(s)=\left(\frac{2\pi is}{e^u-1}\right)^\mu e^{zu+\mu s},\quad u=2\pi i\left(e^s-1\right),
\end{equation}
and $\calC_+$ is the image of $\calL_+$.  $\calC_+$  is a contour  that encircles the origin in the clockwise fashion. 

To obtain an asymptotic expansion we apply Watson's lemma for loop integrals, see \cite[p.~120]{Olver:1997:ASF}. We expand
\begin{equation}\label{gsexp}
g(s)=\sum_{k=0}^\infty g_k s^k,
\end{equation}
substitute this in \eqref{Ip}, and interchange summation and integration. This gives
\begin{equation}\label{Ipas1}
I_+\sim n! \,\frac{e^{2\pi iz}}{(2\pi i)^n}\,\sum_{k=0}^\infty  g_k F_k,
\end{equation}
where
\begin{equation}\label{Fk}
F_k=\frac{1}{2\pi i}\int_{\calC_+} s^{k-\mu}e^{-ns}\,ds,
\end{equation}
with $\calC_+$  extended  to $+\infty$. That is,  we start  the integration along the contour $\calC_+$ at $s=+\infty$, with $\phase\,s=2\pi$, turn around the origin in the clock-wise direction, and return to $+\infty$ with $\phase\,s=0$.

To evaluate  the integrals we turn the path by writing $s=e^{\pi i}t$, and use the representation of the reciprocal gamma function in terms of the Hankel contour; see \cite[p.~48]{Temme:1996:SFI}. The result is
\begin{equation}\label{Fkgam}
F_k=n^{\mu-k-1}e^{\pi i \mu}\frac{(-1)^k}{\Gamma(\mu-k)}=n^{\mu-k-1}e^{\pi i \mu}\frac{(1-\mu)_k}{\Gamma(\mu)},
\end{equation}
where $(a)_n$ is the shifted factorial, or Pochhammer's symbol, defined by
\begin{equation}\label{poch}
(a)_n=a\cdot(a+1)\cdots(a+n-1)=\frac{\Gamma(a+n)}
{\Gamma(a)}, \quad n=0,1,2,\ldots\,.
\end{equation}
This gives the expansion
\begin{equation}\label{Ipas2}
I_+\sim \frac{n! \,n^{\mu-1}}{(2\pi )^n\,\Gamma(\mu)}\,e^{i\chi}\,\sum_{k=0}^\infty  \frac{(1-\mu)_k g_k}{n^k},  
\end{equation}
where
\begin{equation}\label{chizeta}
\chi=2\zeta-\tfrac12n\pi,\quad \zeta=(z+\tfrac12\mu)\pi.
\end{equation}

The contribution $I_-$ can be obtained in a similar way. However, it is the complex conjugate of $I_+$ (not considering $z$ and $\mu$ as complex numbers). The result $I_++I_-$ can be obtained by taking twice the real part of $I_+$. We write $g_k=g_k^{(r)}+ig_k^{(i)}$ (with $g_k^{(r)},g_k^{(i)}$ real when $z$ and $\mu$ are real), and obtain
\begin{equation}\label{Beras}
B_n^\mu(z)\sim \frac{2\,n! \,n^{\mu-1}}{(2\pi )^n\,\Gamma(\mu)}\left[
\cos\chi\,\sum_{k=0}^\infty  \frac{(1-\mu)_k g_k^{(r)}}{n^k}-
\sin\chi\,\sum_{k=0}^\infty  \frac{(1-\mu)_k g_k^{(i)}}{n^k}\right],
\end{equation}
as $n\to\infty$, with $z$ and $\mu$ fixed complex numbers ($\mu\notin\ZZ$).

The first few coefficients $g_k^{(r)},g_k^{(i)}$ are
\renewcommand{\arraystretch}{1.5}
\begin{equation}\label{gkri}
\begin{array}{ll}
g_0^{(r)}=1,  &g_0^{(i)}=0,\\
g_1^{(r)}=\tfrac12\mu,  &g_1^{(i)}=2\zeta,\\
g_2^{(r)}=\tfrac1{24}(3\mu^2+(4\pi^2-1)\mu-48\zeta^2),  &g_2^{(i)}=(1+\mu)\zeta,\\
g_3^{(r)}=\tfrac1{48}(\mu^3+(4\pi^2-1)\mu^2+8(\pi^2-6\zeta^2)\mu-96\zeta^2),  &\\
g_3^{(i)}=\tfrac1{12}\zeta(3\mu^2+(4\pi^2+5)\mu-16\zeta^2+4).&\\
\end{array}
\end{equation}
\renewcommand{\arraystretch}{1.0}

The first-order approximation reads
\begin{equation}\label{Berasf}
B_n^\mu(z)= \frac{2\,n! \,n^{\mu-1}}{(2\pi )^n\,\Gamma(\mu)}\left[
\cos\pi(2z+\mu-\tfrac12n)+\bigO(1/n)\right], \quad n\to\infty.
\end{equation}
N\"orlund \cite[p.~39]{Norlund:1961:SVA} describes the same method of this section and only gives the first-order approximation.

\subsubsection{An alternative expansion}\label{subsec:Ber4}
As observed in the previous method, the main contributions to \eqref{Berint} comes from the singular points of the integrand at $\pm2\pi i$. In this section we expand part of the integrand of  \eqref{Berint} in a two-point Taylor expansion.  In this way a simpler asymptotic representation can be obtained. For more details on this topic we refer to \cite{Lopez:2002:TPT,Lopez:2004:MPT} and for the evaluation of coefficients of such expansions to \cite{Vidunas:2002:SEC}. We write
\begin{equation}\label{Ber41}
f(w)=2^{-3\mu}\pi^{-2\mu}\left(w^2+4\pi^2\right)^\mu\left(\frac{w}{e^w-1}\right)^\mu e^{wz}
\end{equation}
and expand
\begin{equation}\label{Ber42}
f(w)
=\sum_{k=0}^\infty\left(\alpha_k+w\beta_k\right)\left(w^2+4\pi^2\right)^k.
\end{equation}
The function $f(w)$ is analytic inside the disk $|w|<4\pi$ and the series converges in the same domain. The coefficients 
$\alpha_0$ and $\beta_0$ can be found by substituting $w=\pm 2\pi i$. This gives
\renewcommand{\arraystretch}{1.5}
\begin{equation}\label{Ber43}
\begin{array}{l}
\dsp{\alpha_0=\frac{f(2\pi i)+f(-2\pi i)}{2}=\cos2\zeta}, \\ 
\dsp{ \beta_0=\frac{f(2\pi i)-f(-2\pi i)}{4\pi i}=\frac{1}{2\pi}\sin2\zeta, }
\end{array}
\renewcommand{\arraystretch}{1.0}
\end{equation}
where $\zeta$ is defined in \eqref{chizeta}.

The next coefficients can be obtained by writing $f_0(w)=f(w)$ and
\renewcommand{\arraystretch}{1.5}
\begin{equation}\label{Ber44}
\begin{array}{lll}
f_{j+1}(w)&=&\dsp{\frac{f_j(w)-(\alpha_j+w\beta_j)}{w^2+4\pi^2}} \\ 
&=&\dsp{\sum_{k=j+1}^\infty\left(\alpha_k+w\beta_k\right)\left(w^2+4\pi^2\right)^{k-j-1},}
\end{array}
\renewcommand{\arraystretch}{1.0}
\end{equation}
$ j=0,1,2,\ldots$, and by taking the limits when $w\to\pm2\pi i$. We have
\renewcommand{\arraystretch}{1.5}
\begin{equation}\label{Ber45}
\begin{array}{l}
\dsp{\alpha_{j+1}=\frac{f_j^\prime(2\pi i)-f_j^\prime(-2\pi i)}{8\pi i}}, \\ 
\dsp{ \beta_{j+1}=-\frac{f_j^\prime(2\pi i)+f_j^\prime(-2\pi i)-2\beta_j}{16\pi^2}}.
\end{array}
\renewcommand{\arraystretch}{1.0}
\end{equation}
This gives
{\small
\renewcommand{\arraystretch}{1.75}
\begin{equation}\label{Ber46}
\begin{array}{l}
\dsp{ \alpha_1=-\frac{1}{16\pi^2}[3\mu\cos2\zeta+2\pi\eta\sin2\zeta],}\\
\dsp{ \beta_1=\frac{1}{32\pi^3}[2\pi\eta\cos2\zeta+(2-3\mu)\sin2\zeta],}\\
\dsp{\alpha_2=\frac{1}{1536\pi^4}[(-12\pi^2\eta^2+4\mu\pi^2-33\mu+27\mu^2)\cos2\zeta+12\pi\eta(3\mu-1)\sin2\zeta]}, \\ 
\dsp{ \beta_2=\frac{1}{3072\pi^5}[-36\pi\eta(\mu-1)\cos2\zeta+(36-69\mu+27\mu^2+4\mu\pi^2-12\pi^2\eta^2)\sin2\zeta],}
\end{array}
\renewcommand{\arraystretch}{1.0}
\end{equation}
}
where $\eta=\mu-2z$.

Substituting the expansion in \eqref{Ber42} into \eqref{Berint} we obtain
\begin{equation}\label{Ber47}
B_n^{\mu}(z)=n!\, 2^{3\mu}\pi^{2\mu}\sum_{k=0}^\infty \left[\alpha_k \Phi_k^{(n)}+\beta_k \Phi_k^{(n-1)}\right],
\end{equation}
where
\begin{equation}\label{Ber48}
\Phi_k^{(n)}=\frac{1}{2\pi i} \int_{\cal C} \left(w^2+4\pi^2\right)^{k-\mu}\frac{dw}{w^{n+1}}.
\end{equation}
We have $\Phi_k^{(2n+1)}=0$ and
\begin{equation}\label{Ber49}
 \Phi_k^{(2n)}=(2\pi)^{2k-2\mu-2n}\,{k-\mu\choose  n}=(-1)^n(2\pi)^{2k-2\mu-2n}\frac{(\mu-k)_n}{n!}.
\end{equation}
Hence,
\renewcommand{\arraystretch}{1.5}
\begin{equation}\label{Ber410}
\begin{array}{l}
\dsp{B_{2n}^{\mu}(z)=(2n)! \,2^{3\mu}\pi^{2\mu}\sum_{k=0}^\infty \alpha_k \Phi_k^{(2n)},} \\
\dsp{B_{2n+1}^{\mu}(z)=(2n+1)!2^{3\mu}\pi^{2\mu}\,\sum_{k=0}^\infty \beta_k \Phi_k^{(2n)}.}
\end{array}
\renewcommand{\arraystretch}{1.0}
\end{equation}
These convergent expansions have an asymptotic character for large $n$. This follows from (see \eqref{poch})
\renewcommand{\arraystretch}{1.75}
\begin{equation}\label{Ber411}
\begin{array}{@{}r@{\;}c@{\;}l@{}}
 \dsp{\frac{ \Phi_{k+1}^{(2n)}}{ \Phi_k^{(2n)}}}
&=& \dsp{4\pi^2\frac{(\mu-k-1)_n}{(\mu-k)_n}}\\
&=& \dsp{4\pi^2\frac{\Gamma(\mu-k-1+n)}{\Gamma(\mu-k-1)}\,\frac{\Gamma(\mu-k)}{\Gamma(\mu-k+n)}}\\
&=& \dsp{4\pi^2\frac{\mu-k-1}{\mu-k+n-1}=\bigO\left(n^{-1}\right), \quad n\to\infty.}
\end{array}
\renewcommand{\arraystretch}{1.0}
\end{equation}

We compare the first term approximations given in \eqref{Berasf} and those from \eqref{Ber410}.
From \eqref{Berasf}  we obtain
\begin{equation}\label{Ber414}
B_{2n}^\mu(z)\sim (-1)^n\frac{(2n)! \,2^\mu n^{\mu-1}}{(2\pi )^{2n}\,\Gamma(\mu)}
\cos\pi(2z+\mu)+\ldots, 
\end{equation}
and from  \eqref{Ber410}
\begin{equation}\label{Ber415}
B_{2n}^\mu(z)=(-1)^n\frac{(2n)! \,2^\mu }{(2\pi )^{2n}\,\Gamma(\mu)}\frac{\Gamma(n+\mu)}{n!}
\cos\pi(2z+\mu)+\ldots.
\end{equation}
Because $\Gamma(n+\mu)/n!\sim n^{\mu-1}$ as $n\to\infty$, we see that the first approximations give the same asymptotic estimates, and they are exactly the same when $\mu=1$.

\paragraph{Integer values of \protectbold{\mu}}\label{par:Bintmu}
Comparing the expansions in \eqref{Beras} and \eqref{Ber410}, we observe that those in \eqref{Ber410} do not vanish when $\mu=0,-1,-2,\ldots$, whereas the expansion in \eqref{Beras} does. We have when $\mu=m$ (integer)
\renewcommand{\arraystretch}{1.5}
\begin{equation}\label{Ber412}
\Phi_k^{(2n)}=
\left\{
\begin{array}{ll}
\dsp{(2\pi)^{2k-2m-2n}\,{k-m\choose  n}}, \quad &k\ge n+m,\\
0, \quad & k<n+m.
\end{array}
\right.
\renewcommand{\arraystretch}{1.0}
\end{equation}
Hence, the summation in \eqref{Ber410} starts with $k=n+m$. The scale $\{\Phi_k^{2n}\}$ loses its asymptotic property, because now
\begin{equation}\label{Ber413}
\frac{ \Phi_{k+1}^{(2n)}}{ \Phi_k^{(2n)}}=4\pi^2\frac{n+\ell+1}{\ell+1}=\bigO\left(n\right), \quad n\to\infty,
\end{equation}
where $k=n+m+\ell$, and a possible asymptotic character of the series in \eqref{Ber410} has  to be furnished by the coefficients $\alpha_k, \beta_k$, which depend on $n$ when $k\ge n+m$

Because $n$ is assumed to be large, and the coefficients $\alpha_k, \beta_k$ in  \eqref{Ber410} become quite complicated when $k\ge n+m$, these expansions are of no use when $\mu$ is an integer.

When we replace the expansion in \eqref{Ber42} with
\begin{equation}\label{Ber42a}
f(w)
=\sum_{k=0}^\infty\left(\wt\alpha_k+w\wt\beta_k\right)\left(\frac{w^2+4\pi^2}{w^2}\right)^k
\end{equation}
we obtain
\renewcommand{\arraystretch}{1.5}
\begin{equation}\label{Ber410a}
\begin{array}{l}
\dsp{B_{2n}^{\mu}(z)\sim(2n)! \,2^{3\mu}\pi^{2\mu}\sum_{k=0}^\infty \wt\alpha_k \wt\Phi_k^{(2n)},} \\
\dsp{B_{2n+1}^{\mu}(z)\sim(2n+1)!2^{3\mu}\pi^{2\mu}\,\sum_{k=0}^\infty\wt\beta_k \wt\Phi_k^{(2n)},}
\end{array}
\renewcommand{\arraystretch}{1.0}
\end{equation}
where $\wt\alpha_k$ and $\wt\beta_k$ can be obtained from a similar scheme as in \eqref{Ber44}-\eqref{Ber45}. The functions $\Phi_k^{(2n)}$ are given by
\begin{equation}\label{Ber49a}
\wt \Phi_k^{(2n)}=\frac{(-1)^{n+k}}{(2\pi)^{2\mu+2n}}\frac{(\mu-k)_{n+k}}{(n+k)!}
=\frac{(-1)^{n+k}}{(2\pi)^{2\mu+2n}}\frac{\Gamma(\mu+n)}{\Gamma(\mu-k)\,(n+k)!}.
\end{equation}
When $\mu=m$ (integer) these functions vanish if $k-m=0,1,2,\ldots$, which is more useful that in the earlier choice \eqref{Ber42}. When $m<0$ all terms vanish, when $m>0$ the series have a finite number of terms. 

With the expansion in \eqref{Ber42a}, which converges in certain neighborhoods of the points $w=\pm2\pi i$ (and not in a domain that contains any allowed deformation of  the curve $\calC$ in \eqref{Berint}), the expansions in \eqref{Ber410a} do not converge, but they have an asymptotic character for large $n$. 

As an example, when $m=1$ the expansions in \eqref{Ber410a} have just one term $(k=0)$. In this case
$\wt\Phi_0^{(2n)}=(-1)^n/(2\pi)^{2n+2}$ and $\wt\alpha_0=\alpha_0$, $\wt\beta_0=\beta_0$ (see \eqref{Ber43}). These approximations correspond exactly to the first terms in the expansions in \eqref{Bersereven} and \eqref{Berserodd}.

\section{The generalized Euler polynomials}\label{sec:Eul}
We can use the same methods as for the Bernoulli polynomials, and, therefore, we give less details. Again, three different cases arise, depending on $\mu=0,-1,-2,\ldots$, $\mu=1,2,3,\dots$, and $\mu$ otherwise, real or complex. In the first and third case we use the Cauchy integral
\begin{equation}\label{Eulint}
E_n^{\mu}(z)=\frac{n!}{ 2\pi i} \int_{\calC}\frac{2^{\mu}e^{wz}}{(e^w+1)^{\mu}}
\frac{dw}{ w^{n+1}},
\end{equation}
where $\calC$ is a circle around the origin, with radius less than $\pi$. This follows from \eqref{Euldef}.

\subsection{Asymptotic form when  \protectbold{\mu=0,-1,-2,\ldots}}\label{subsec:Eul1}
We proceed as in \S\ref{subsec:Ber1} and write $\mu=-m$ $(m=0,1,2,\ldots)$. We expand $E_n^{\mu}(z)$ in terms of a finite sum. We have
\begin{equation}\label{Eulsum}
E_n^{-m}(z)=2^{-m} \sum_{r=0}^m {m\choose r} (z+r)^{n}.
\end{equation}
For any given $m\in \NN$ and complex $z$ only the term or terms with the largest values of $|z+r|$ will give a large contribution to the sum in \eqref{Eulsum}, the other terms being exponentially small in comparison. We conclude, that \eqref{Eulsum} gives the asymptotic form when $n\to\infty$, when $m$ and $z$ are fixed. In particular, when $z>0$, the term with index $r=m$ is maximal, and we have
\begin{equation}\label{Eulsum0}
E_n^{-m}(z)=2^{-m}  (z+m)^{n}\left[1+\bigO\left(\frac{z+m-1}{z+m}\right)^{n}\right].
\end{equation}
The error term can also be estimated by $\bigO(\exp(-n/(z+m)))$, which is indeed exponentially small compared with unity.

For general complex $z=x+iy$ and $x>-m/2$ the term with index $r=m$ again is maximal and the same estimate as in \eqref{Eulsum0} is valid. When $x=-m/2$ the terms with $r=0$ and $r=m$ give the maximal contributions, and we have
 \begin{equation}\label{Eulsum1}
E_n^{-m}(z)\sim2^{-m}  \left[(-\tfrac12m+iy)^{n}+(\tfrac12m+iy)^{n} \right].
\end{equation}
When $x<-m/2$ the term with index $r=0$ is maximal, and we have
\begin{equation}\label{Eulsum2}
E_n^{-m}(z)=2^{-m}z^{n}\left[1+\bigO\left(\frac{z+1}{z}\right)^{n}\right].
\end{equation}

As explained at the end of \S\ref{subsec:Ber1} these estimates can also be derived by using the saddle point method.

\subsection{Asymptotic form when  \protectbold{\mu= 1, 2,3,\ldots}}\label{subsec:Eul2}
We write $\mu=m$. For $m=1$ we have
\renewcommand{\arraystretch}{1.75}
\begin{equation}\label{Eul1ser}
\begin{array}{@{}r@{\;}c@{\;}l@{}}
E_n^{1}(z)=E_n(z)
&=& \dsp{2\,n!\sum_{{k=-\infty}}^\infty\frac{e^{(2k+1)\pi iz}}{((2k+1)\pi i)^{n+1}}}\\
&=& \dsp{4\,n!\sum_{{k=0}}^\infty\frac{\sin((2k+1)\pi z-\frac12\pi n)}{((2k+1)\pi)^{n+1}}},
\end{array}
\renewcommand{\arraystretch}{1.0}
\end{equation}
where $z\in(0,1)$ if $n=0$ and $z\in[0,1]$ if $n>0$. This expansion follows from   \eqref{Eulint} as in \S\ref{subsec:Ber2}.

In the second series in \eqref{Eul1ser} only the term  with $k=0$ is relevant for the asymptotic behavior, and we obtain for fixed complex $z$
\begin{equation}\label{Eul1as}
E_n^{1}(z)=\frac{4\,n!}{\pi^{n+1}}\left[\sin\left(\pi z-\tfrac12\pi n\right) +\bigO\left(3^{-n}\right)\right],\quad n\to\infty.
\end{equation}

For general fixed real or complex $z$ the series in  \eqref{Eul1ser}  can be viewed as asymptotic expansion for large $n$, as easily follows from the ratio test.

For general $\mu=m=1,2,3,\ldots$ a similar expansion as in \eqref{Eul1ser} can be given. In that case the poles of the integrand in  \eqref{Eulint} are of higher order. We can write
\begin{equation}\label{Eulmser}
E_n^{m}(z)=2\,n!\sum_{{k=0}}^\infty\epsilon_k^m(n,z)\frac{e^{(2k+1)\pi iz}}{((2k+1)\pi i)^{n+1}},
\end{equation}
where $\epsilon_k^1(n,z)=1, \forall k$.  

To obtain  $\epsilon_k^m(n,z)$ for  $m> 1$ we compute the residues of the poles at $(2k+1)\pi i$ of order $m$ of the integrand in \eqref{Eulint}. For this we compute the coefficient $d_{m-1}$ in the expansion
\begin{equation}\label{Eulres1}
\frac{(w-(2k+1)\pi i)^m e^{zw}}{(e^w+1)^mw^{n+1}}=\sum_{r=0}^\infty d_r (w-(2k+1)\pi i)^r.
\end{equation}
We substitute  $w=s+(2k+1)\pi i$ and write the expansion as
\begin{equation}\label{Eulres2}
(-1)^me^{z(2k+1)\pi i}\frac{s^m e^{zs}}{(e^s-1)^m(s+(2k+1)\pi i)^{n+1}}=\sum_{r=0}^\infty d_r s^r.
\end{equation}
We use \eqref{Berdef} and conclude that 
\begin{equation}\label{Eulres3}
d_{m-1}=\frac{(-1)^me^{z(2k+1)\pi i}}{((2k+1)\pi i)^{n+1}}\sum_{\nu=0}^{m-1} \frac{B_\nu^m(z)}{\nu!}{-n-1\choose m-1-\nu}((2k+1)\pi i)^{\nu+1-m}.
\end{equation}
It follows that 
\begin{equation}\label{Eulres4}
\epsilon_k^m(n,z)=(-1)^{m-1}2^{m-1}\sum_{\nu=0}^{m-1} \frac{B_\nu^m(z)}{\nu!}{-n-1\choose m-1-\nu}((2k+1)\pi i)^{\nu+1-m}, 
\end{equation}
which we write in the form
\renewcommand{\arraystretch}{1.5}
\begin{equation}\label{Eulres5}
\begin{array}{l}
\dsp{\epsilon_k^m(n,z)=\frac{2^{m-1}}{((2k+1)\pi i)^{m-1}}{n+m-1\choose m-1}}\,\times \\
\quad\quad
\dsp{\sum_{\nu=0}^{m-1}B_\nu^m(z){m-1\choose \nu}\frac{(n+m-\nu-1)!}{(n+m-1)!}(-(2k+1)\pi i)^{\nu}}.
\end{array}
\end{equation}
\renewcommand{\arraystretch}{1.0}

For large $n$ the main term occurs for $\nu=0$, giving
\begin{equation}\label{Eulres6}
\epsilon_k^m(n,z)=
\frac{2^{m-1}n^{m-1}}{(m-1)!\,((2k+1)\pi i)^{m-1}}\left[1+\bigO\left(n^{-1}\right)\right],
\end{equation}
and in \eqref{Eulmser} the terms with $k=0, -1$ give the main terms, and we obtain for fixed $m$ and complex $z$ (cf. \eqref{Eul1as})
\renewcommand{\arraystretch}{1.5}
\begin{equation}\label{Eulmas1}
\begin{array}{ll}
\dsp{E_n^{m}(z)=\frac{2^{m+1}n!}{\pi^{n+m}}{n+m-1\choose m-1} \,\times} \\
\quad\quad
\dsp{\left[
\sum_{\nu=0}^{m-1} B_\nu^m(z)
{m-1\choose \nu}\frac{(n+m-\nu-1)!}{(n+m-1)!}\pi^{\nu}\sin\tau+\bigO\left(3^{-n}\right)\right]},
\end{array}
\renewcommand{\arraystretch}{1.0}
\end{equation}
as $n\to\infty$, where $\tau=(z-\tfrac12n-\tfrac12(m-1)-\tfrac12\nu)\pi$.

\subsection{Asymptotic form for general complex  \protectbold{\mu}}\label{subsec:Eul3}
The analysis is as in \S\ref{subsec:Ber3}. We use a contour for the integral \eqref{Eulint} as in Figure~\ref{Ber.fig1}, now with loops around the branch points $\pm\pi i$, and with radius of the large circle smaller  than $3\pi$. We denote the integrals around the loops by $I_\pm$. After the substitution $w=\pi i\exp(s)$ we obtain for the upper loop
\begin{equation}\label{IpEul}
I_+=\frac{2^\mu n!}{2\pi i}\frac{e^{\pi iz-\mu\pi i}}{(\pi i)^{n+\mu}}\int_{\calC_+} h(s) s^{-\mu}e^{-ns}\,ds,
\end{equation}
where
\begin{equation}\label{hs}
h(s)=e^{zu}\left(\frac{\pi i s}{e^u-1}\right)^\mu,\quad u=\pi i \left(e^s-1\right).
\end{equation}
We expand
$h(s)=\sum_{k=0}^\infty h_k s^k $ and interchange summation and integration in \eqref{IpEul}.
By using \eqref{Fk} and \eqref{Fkgam} we obtain the result
\begin{equation}\label{Eulas}
E_n^\mu(z)\sim\frac{2^{\mu+1}\,n! \,n^{\mu-1}}{\pi ^{n+\mu}\,\Gamma(\mu)}\left[
\cos\chi\,\sum_{k=0}^\infty  \frac{(1-\mu)_k h_k^{(r)}}{n^k}-
\sin\chi\,\sum_{k=0}^\infty  \frac{(1-\mu)_k h_k^{(i)}}{n^k}\right],
\end{equation}
as $n\to\infty$, with $z$ and $\mu$ fixed complex numbers ($\mu\notin\ZZ$),
where
\begin{equation}\label{chizeta2}
\chi=\zeta-\tfrac12n\pi,\quad \zeta=(z-\tfrac12\mu)\pi.
\end{equation}
The first few coefficients $h_k^{(r)}, h_k^{(i)}$ are
\renewcommand{\arraystretch}{1.5}
\begin{equation}\label{hkri}
\begin{array}{ll}
h_0^{(r)}=1,  %\quad  
&h_0^{(i)}=0,\\
%\quad 
h_1^{(r)}=-\tfrac12\mu, % \quad 
&h_1^{(i)}=\zeta,\\
h_2^{(r)}=\tfrac1{24}(3(1-2\pi^2)\mu^2+(13\pi^2-12\zeta\pi-1)\mu-12\zeta^2),   %\\
&h_2^{(i)}=\tfrac12(1-\mu)\zeta, \\
h_3^{(r)}=\tfrac{1}{48}z(-\mu^3 +(1-\pi^2)\mu^2+2(\pi^2+6\zeta^2)\mu-24\zeta^2), &\\
h_3^{(i)}=\tfrac1{24}\zeta(3\mu^2+(\pi^2-7)\mu-4\zeta^2+4).&\\
\end{array}
\end{equation}
\renewcommand{\arraystretch}{1.0}

The first-order approximation reads
\begin{equation}\label{Eulasf}
E_n^\mu(z)=  \frac{2^{\mu+1}\,n! \,n^{\mu-1}}{\pi ^{n+\mu}\,\Gamma(\mu)}\left[
\cos\pi(z-\tfrac12\mu-\tfrac12n)+\bigO(1/n)\right],\quad n\to\infty.
\end{equation}

\subsubsection{An alternative expansion}\label{subsec:Eul4}
We repeat the steps of  \S\ref{subsec:Ber4}. We write
\begin{equation}\label{Eul41}
g(w)=\left(\frac{w^2+\pi^2}{2\pi}\right)^\mu\left(\frac{1}{e^w+1}\right)^\mu e^{wz}
\end{equation}
and expand
\begin{equation}\label{Eul42}
g(w)
=\sum_{k=0}^\infty\left(\gamma_k+w\delta_k\right)\left(w^2+\pi^2\right)^k.
\end{equation}
We have
\renewcommand{\arraystretch}{1.5}
\begin{equation}\label{Eul43}
\begin{array}{l}
\dsp{\gamma_0=\frac{g(\pi i)+g(-\pi i)}{2}=\cos \zeta}, \\ 
\dsp{ \delta_0=\frac{g(\pi i)-g(-\pi i)}{2\pi i}=\frac{1}{\pi}\sin \zeta. }
\end{array}
\renewcommand{\arraystretch}{1.0}
\end{equation}
where $\zeta=(z-\frac12\mu)\pi$.

The next coefficients follow from writing $g_0(w)=g(w)$ and 
\renewcommand{\arraystretch}{1.5}
\begin{equation}\label{Eul44}
\begin{array}{lll}
g_{j+1}(w)&=&\dsp{\frac{g_j(w)-(\gamma_j+w\delta_j)}{w^2+\pi^2}} \\ 
&=&\dsp{\sum_{k=j+1}^\infty\left(\gamma_k+w\delta_k\right)\left(w^2+\pi^2\right)^{k-j-1}, \quad j=0,1,2,\ldots.}
\end{array}
\renewcommand{\arraystretch}{1.0}
\end{equation}
This gives
\renewcommand{\arraystretch}{1.5}
\begin{equation}\label{Eul45}
\begin{array}{lll}
\dsp{\gamma_{j+1}=\frac{g_j^\prime(\pi i)-g_j^\prime(-\pi i)}{4\pi i}}, \\ 
\dsp{ \delta_{j+1}=-\frac{g_j^\prime(\pi i)+g_j^\prime(-\pi i)-2\delta_j}{4\pi^2}}.
\end{array}
\renewcommand{\arraystretch}{1.0}
\end{equation}
and
{\small
\renewcommand{\arraystretch}{1.75}
\begin{equation}\label{Eul46}
\begin{array}{l}
\dsp{ \gamma_1=-\frac{1}{4\pi^2}[\mu\cos\zeta+\pi\eta\sin  \zeta],}\\
\dsp{ \delta_1=\frac{1}{4\pi^3}[\pi\eta\cos\zeta+(2-\mu)\sin\zeta],}\\
\dsp{\gamma_2=\frac{1}{96\pi^4}[(-9\mu-3\pi^2\eta^2+\pi^2\mu+3\mu^2)\cos\zeta+ 6\pi\eta(\mu-1)\sin\zeta]}, \\ 
\dsp{ \delta_2=\frac{1}{96\pi^5}[6\pi\eta(3-\mu)\cos\zeta+(36-21\mu+3\mu^2+\pi^2\mu-3\pi^2\eta^2)\sin\zeta],}
\end{array}
\renewcommand{\arraystretch}{1.0}
\end{equation}
}
where $\eta=\mu-2z$.

Substituting the expansion in \eqref{Eul42} into \eqref{Eulint} we obtain
\begin{equation}\label{Eul47}
E_n^{\mu}(z)=(4\pi)^\mu n!\sum_{k=0}^\infty \left[\gamma_k \Psi_k^{(n)}+\delta_k \Psi_k^{(n-1)}\right],
\end{equation}
where
\begin{equation}\label{Eul48}
\Psi_k^{(n)}=\frac{1}{2\pi i} \int_{\cal C} \left(w^2+\pi^2\right)^{k-\mu}\frac{dw}{w^{n+1}}.
\end{equation}
We have $\Psi_k^{(2n+1)}=0$ and
\begin{equation}\label{Eul49}
 \Psi_k^{(2n)}=\pi^{2k-2\mu-2n}\,{k-\mu\choose  n}=(-1)^n\pi^{2k-2\mu-2n}\frac{(\mu-k)_n}{n!}.
\end{equation}
Hence,
\renewcommand{\arraystretch}{1.5}
\begin{equation}\label{Eul410}
\begin{array}{l}
\dsp{E_{2n}^{\mu}(z)=(2\pi)^\mu(2n)!\sum_{k=0}^\infty \gamma_k \Psi_k^{(2n)},} \\
\dsp{E_{2n+1}^{\mu}(z)=(2\pi)^\mu(2n+1)!\sum_{k=0}^\infty \delta_k \Psi_k^{(2n)}.}
\end{array}
\renewcommand{\arraystretch}{1.0}
\end{equation}
These convergent expansions have an asymptotic character for large $n$. This follows from
\begin{equation}\label{Eul411}
\frac{ \Psi_{k+1}^{(2n)}}{ \Psi_k^{(2n)}}=\pi^2\frac{\mu-k}{\mu-k-1+n}=\bigO\left(n^{-1}\right), \quad n\to\infty.
\end{equation}

Comparing the first term approximations given in \eqref{Eulasf} and those from \eqref{Eul410}
 we obtain from  \eqref{Eulasf}
\begin{equation}\label{Eul413}
E_{2n}^\mu(z)\sim (-1)^n\frac{(2n)! \,2^{2\mu} n^{\mu-1}}{\pi^{2n+\mu}\,\Gamma(\mu)}
\cos\pi(z-\tfrac12\mu)+\ldots, 
\end{equation}
and from  \eqref{Eul410}
\begin{equation}\label{Eul414}
E_{2n}^\mu(z)=(-1)^n\frac{(2n)! \,2^{2\mu} }{\pi^{2n+\mu}\,\Gamma(\mu)}\frac{\Gamma(n+\mu)}{n!}
\cos\pi(z-\tfrac12\mu)+\ldots.
\end{equation}
and we see that the first approximations give the same asymptotic estimates.

\paragraph{Integer values of \protectbold{\mu}}\label{par:Eintmu}
The expansions in \eqref{Eul410} do not vanish when $\mu$ is a negative integer, as the expansion in \eqref{Eulas} does. We have when $\mu=m$ (integer)
\renewcommand{\arraystretch}{1.5}
\begin{equation}\label{Eul412}
\Psi_k^{(2n)}=
\left\{
\begin{array}{ll}
\dsp{\pi^{2k-2m-2n}\,{k-m\choose  n}}, \quad &k\ge n+m,\\
0, \quad & k<n+m.
\end{array}
\right.
\renewcommand{\arraystretch}{1.0}
\end{equation}
Hence, the summation in \eqref{Eul410} starts with $k=n+m$. 

When we expand
\begin{equation}\label{Eul42a}
g(w)
=\sum_{k=0}^\infty\left(\wt\gamma_k+w\wt\delta_k\right)\left(\frac{w^2+\pi^2}{w^2}\right)^k
\end{equation}
we obtain the expansions
\renewcommand{\arraystretch}{1.5}
\begin{equation}\label{Eul410a}
\begin{array}{l}
\dsp{E_{2n}^{\mu}(z)\sim(2\pi)^\mu(2n)!\sum_{k=0}^\infty \wt\gamma_k \wt\Psi_k^{(2n)},} \\
\dsp{E_{2n+1}^{\mu}(z)\sim(2\pi)^\mu(2n+1)!\sum_{k=0}^\infty \wt\delta_k \wt\Psi_k^{(2n)},}
\end{array}
\renewcommand{\arraystretch}{1.0}
\end{equation}
where $\wt\gamma_k$ and $\wt\delta_k$ can be obtained from a similar scheme as in \eqref{Eul45}. The functions $\wt\Psi_k^{(2n)}$ are given by
\begin{equation}\label{Eul49a}
\wt\Psi_k^{(2n)}=\pi^{-2\mu-2n}\,{k-\mu\choose  n+k}=(-1)^{n+k}\pi^{-2\mu-2n}\frac{(\mu-k)_{n+k}}{(n+k)!}.
\end{equation}
When $\mu=m$ (integer) these functions vanish if $k-m=0,1,2,\ldots$, which is more useful that in the earlier choice \eqref{Eul42}. For example, when $m=1$ the expansions in \eqref{Eul410a} have just one term $(k=0)$. In this case
$\wt\Psi_0^{(2n)}=(-1)^n/\pi^{2n+2}$ and $\wt\gamma_0=\gamma_0$, $\wt\delta_0=\delta_0$ (see \eqref{Eul43}). These approximations correspond exactly to the first term in the second expansion in \eqref{Eul1ser}.
 \medskip

\paragraph{Acknowledgments.}
The authors thank the referee for the constructive remarks.
The {\it Gobierno of Navarra, Res. 07/05/2008} is acknowledged for its
financial support.
JLL acknowledges
financial support from {\emph{Ministerio de Educaci\'on y Ciencia}}, project
MTM2007--63772. 
NMT acknowledges financial support from {\emph{Ministerio de Educaci\'on y Ciencia}},
project MTM2006--09050.

%\newpage

\end{document}